\magnification=1200
\input amssym.def
\input epsf

\def \BZ {{\Bbb Z}}

\def\sqr#1#2{{\vcenter{\vbox{\hrule height.#2pt
     \hbox{\vrule width.#2pt height#1pt \kern#1pt
     \vrule width.#2pt} \hrule height.#2pt}}}}
\def\square{\ ${\mathchoice\sqr34\sqr34\sqr{2.1}3\sqr{1.5}3}$}

\centerline {\bf ON THE $z$-DEGREE OF THE KAUFFMAN POLYNOMIAL}
\centerline {\bf OF A TANGLE DECOMPOSITION}
\bigskip
\centerline {Mark E. Kidwell and Theodore B. Stanford}
\centerline {Mathematics Department}
\centerline {United States Naval Academy}
\centerline {572 Holloway Road}
\centerline {Annapolis, MD\ \ 21402}
\medskip
\centerline {\tt mek@nadn.navy.mil}
\centerline {\tt stanford@nadn.navy.mil}
\bigskip

\bigskip
\noindent
{\bf 0. Introduction.}
\bigskip

In 1987, the elder author produced [2] an upper bound on the
degree of the then-new Brandt-Lickorish-Millett-Ho polynomial
in terms of the crossing number and the length of the longest
bridge in a link diagram.  This result extends immediately to
the $z$-degree of the two-variable Kauffman polynomial
(in any of its forms; we shall use the Dubrovnik version
$D(\lambda, z)$).  

The length of the longest bridge (that is, the longest
consecutive string of overcrossings) in a diagram represents
some measure of how nonalternating the diagram is.  For
a fixed crossing number, the greatest $z$-degree will
occur when the diagram is prime and alternating,
which is when
the longest bridge (and every bridge) has length $1$.
A bridge of length $>1$ contributes to a quicker 
unravelling of the link under the skein relations, and
lowers the $z$-degree.  More precisely, 
it is shown in [2] that
the $z$-degree of the polynomial is less than or equal
to $N-B$, where $N$ is the crossing number and $B$ is
the length of the longest bridge in a given link diagram.
It is furthermore shown that if the link diagram is
prime, reduced, and alternating, then the $z$-degree is
$N-1$.

For a diagram with more than one
long bridge, one would like to find an inequality that
considers more than just the single longest bridge.
For example, it was pointed out by Thistlethwaite~[6]
that if a link diagram is composite, then the longest
bridge in each factor counts toward lowering the $z$
degree of $D$.  On the other hand, there are numerous
examples 
(such as $8_{19}$--$8_{21}$ from the table
in Rolfsen~[5])
of non-alternating
diagrams with $N$ crossings and two or more bridges of
length $2$ and with $z$-degree (of the Kauffman polynomial)
of $N-2$.  Thus it
is necessary to find some method of keeping the bridges
separate from each other while computing the skein
tree.  We will accomplish this by cutting the
link diagram into tangles, and considering the longest
bridge in each tangle.  We obtain an upper bound on the
$z$-degree in terms of the number of crossings in the
diagram, and the lengths of each of these longest
separated bridges. 

The Dubrovnik polynomial of a tangle may be defined as
a linear combination, over an appropriate ring, of simple
tangles.  We bound the $z$-degree of the polynomials
in this linear combination in terms of the crossing number
and length of the longest bridge in each tangle.
Our result for links follows by closing up the tangles.

Yokota~[7] has shown that if $L$ is represented by a reduced
alternating diagram with $n$ crossings, 
then the span of the Kauffman 
polynomial in the other variable ($\lambda$ or $a$) is 
equal to $n$.  Our result concerns not the span of $z$ but
the degree.  (The highest negative power of $z$ that occurs
is always one less than the number of components in the link.)

\medskip
\noindent
{\bf Acknowledgement:} \ The authors would like to thank
Joan Birman for her indelible contributions to their
careers.

\bigskip
\noindent
{\bf 1.  The Bound.}
\bigskip

By a {\it tangle} we mean a planar rectangular 
diagram, with overcrossings and undercrossings labeled
in the usual way, containing any number of closed
``circle''
components and exactly two ``arc'' components.  Two
of the arc endpoints are on the top of the rectangle
and two on the bottom.  Equivalence of tangles is up
to regular isotopy.  A {\it bridge} is a maximal segment
of a tangle containing no undercrossings.  The {\it length}
of a bridge is the number of crossings at which it
overcrosses.

Following Morton and Traczyk [4], we shall work over the
ring $\Lambda^\prime$ generated over $\BZ$ by
$\lambda^{\pm 1}, z$, and $\delta$ with the single relation
$$\lambda^{-1} - \lambda = z (\delta-1) \leqno (*)$$
At the end of our calculations, we will use this relation
to eliminate $\delta$ at the expense of introducing
$z^{-1}$ (but not lowering the $z$-degree).

The module $M_2$ is defined to be the set of all 
$\Lambda^\prime$-linear combinations of tangles modulo
the following local relations, where $T \coprod {\rm unknot}$
means the addition of a single unknotted and unlinked circle
component to $T$, and the meaning of
the rest of the symbols is indicated in Figure~1.

\medskip
\item {(i)}
$T^+ - T^- = z (T^0 - T^\infty)$
\medskip
\item {(ii)}
$T^{\rm right} = \lambda^{-1}T$; $T^{\rm left} = \lambda T$
\medskip
\item {(iii)}
$T \coprod {\rm unknot} = \delta T$

\bigskip
\bigskip
\centerline {\hbox {
\epsfysize = 1.5truecm
\epsffile {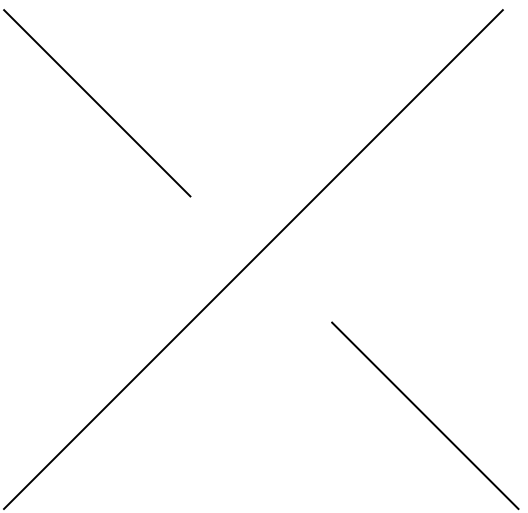} 
\hskip 1truecm
\epsfysize = 1.5truecm
\epsffile {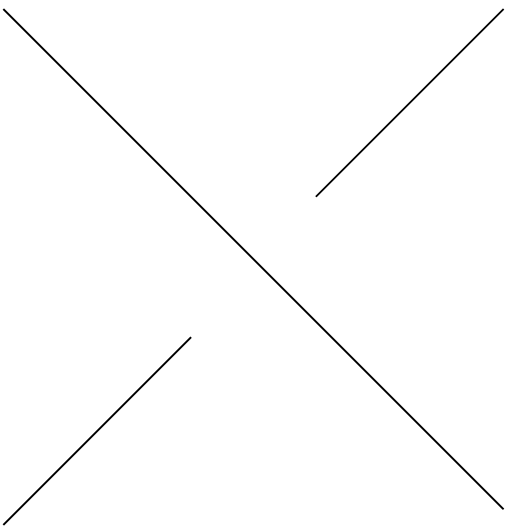} 
\hskip 1truecm 
\epsfysize = 1.5truecm
\epsffile {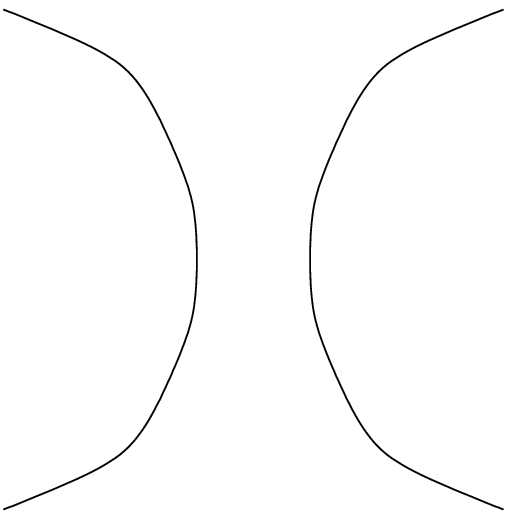} 
\hskip 1truecm
\epsfysize = 1.5truecm
\epsffile {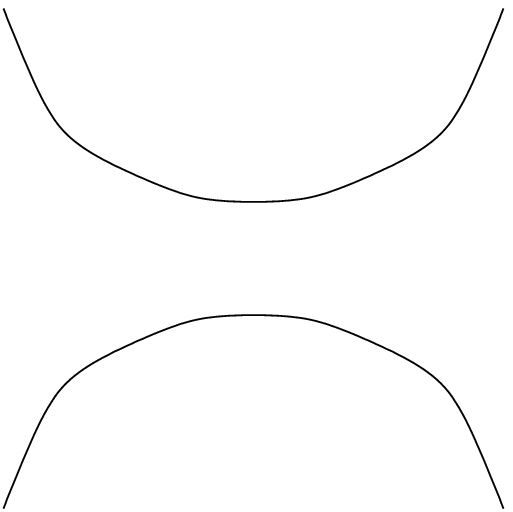}
\hskip 1truecm
\epsfysize = 1.5truecm
\epsffile {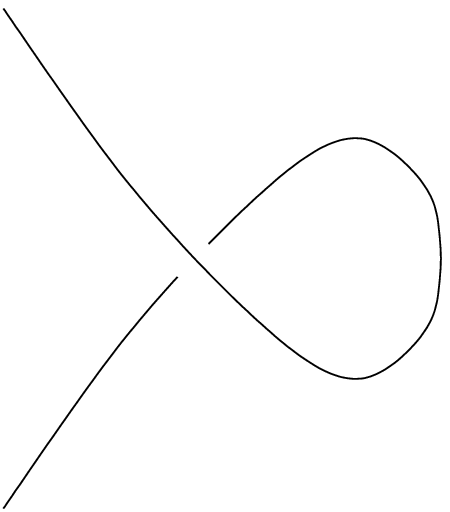} 
\hskip 1truecm
\epsfysize = 1.5truecm
\epsffile {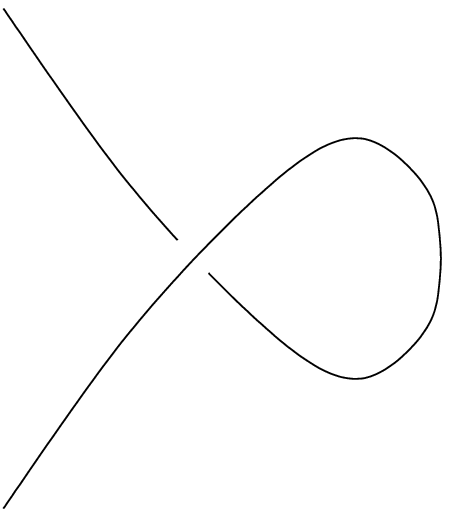}}}

\bigskip
\hbox {
\hskip 2truecm $T^+$
\hskip 1.7truecm $T^-$
\hskip 1.7truecm $T^0$
\hskip 1.7truecm $T^\infty$
\hskip 1.7truecm $T^{\rm left}$
\hskip 1.2truecm $T^{\rm right}$
\hfil}

\bigskip
\centerline {Figure 1}
\bigskip

Morton and Traczyk [4] prove that the tangles called
$P,Q$ and,$R_1$ in Figure~2 form a free $\Lambda^\prime$-basis
for $M_2$.  We find it easier to control the $z$-degree
by working with $P,Q,R_1$, and $R_2$.  
For $f \in \Lambda^\prime$, we understand the $z$-degree
of $f$ to be the minimum $z$-degree over all polynomials
in the equivalence class of $f$.  We shall also sometimes
refer to the $z$-degree of the Kauffman polynomial of a link
or diagram as simply the $z$-degree of the link or diagram.

\bigskip
\centerline {\hbox {
\epsfysize = 2.5truecm
\epsffile {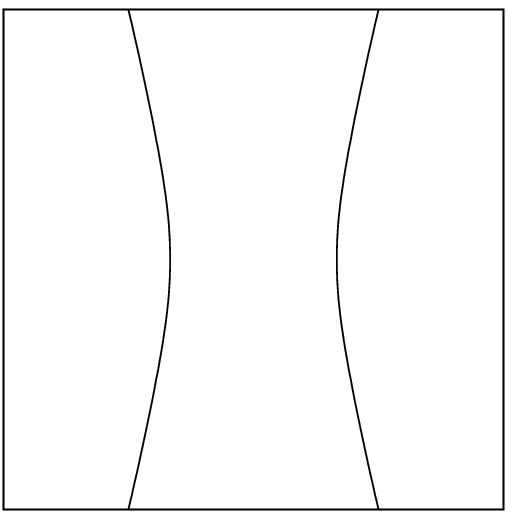} 
\hskip 2truecm
\epsfysize = 2.5truecm
\epsffile {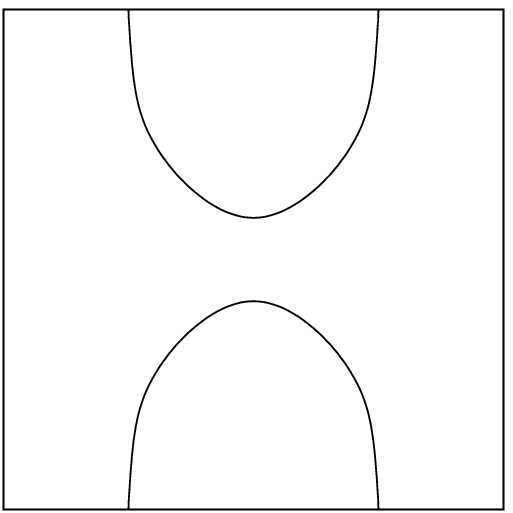} 
\hskip 2truecm
\epsfysize = 2.5truecm
\epsffile {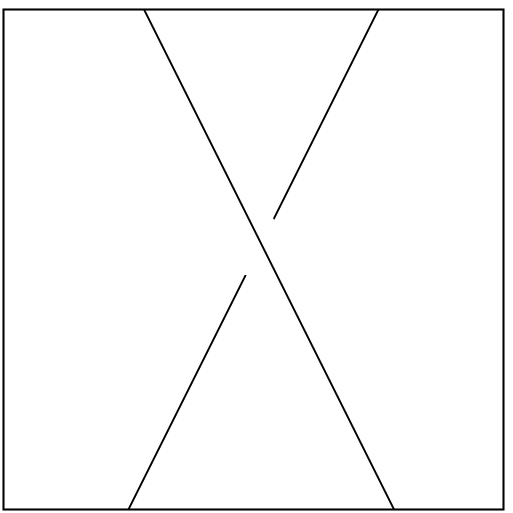} 
\hskip 2truecm
\epsfysize = 2.5truecm
\epsffile {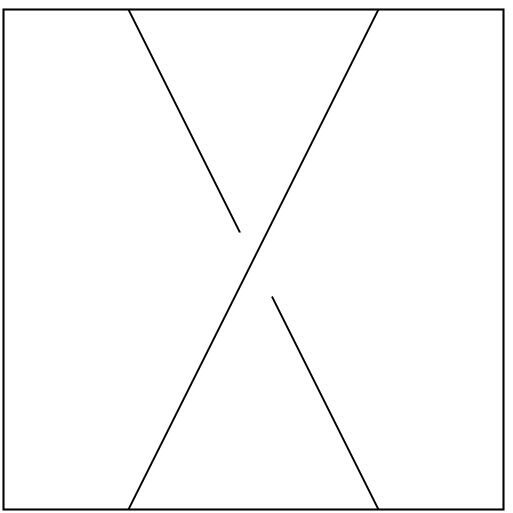}}}
\bigskip

\hbox {
\hskip 1.2truecm $P$
\hskip 4truecm $Q$
\hskip 4truecm $R_1$
\hskip 3.8truecm $R_2$
\hfil}

\bigskip
\centerline {Figure 2}
\bigskip

As in [2], there are a number of bad situations that
must be dealt with before we can proceed to the main
argument.  Call a bridge $b$ in a tangle $T$ {\it improper}
if any of the following is true.  Examples are shown in
Figure~3.

\smallskip
\item {(a)}
The bridge $b$ consists of a full circle with length
$B>0$.
\smallskip
\item {(b)}
One or both of the crossings at which the bridge ends also
has $b$ as an overcrossing.
\smallskip
\item {(c)}
The bridge $b$ starts and ends at the same crossing and 
has length $B>1$.
\smallskip
\item {(d)}
The bridge begins and ends at an endpoint of the tangle
and has length $B>1$.

\bigskip
\bigskip
\centerline {\hbox {
\epsfxsize = 5truecm
\epsffile {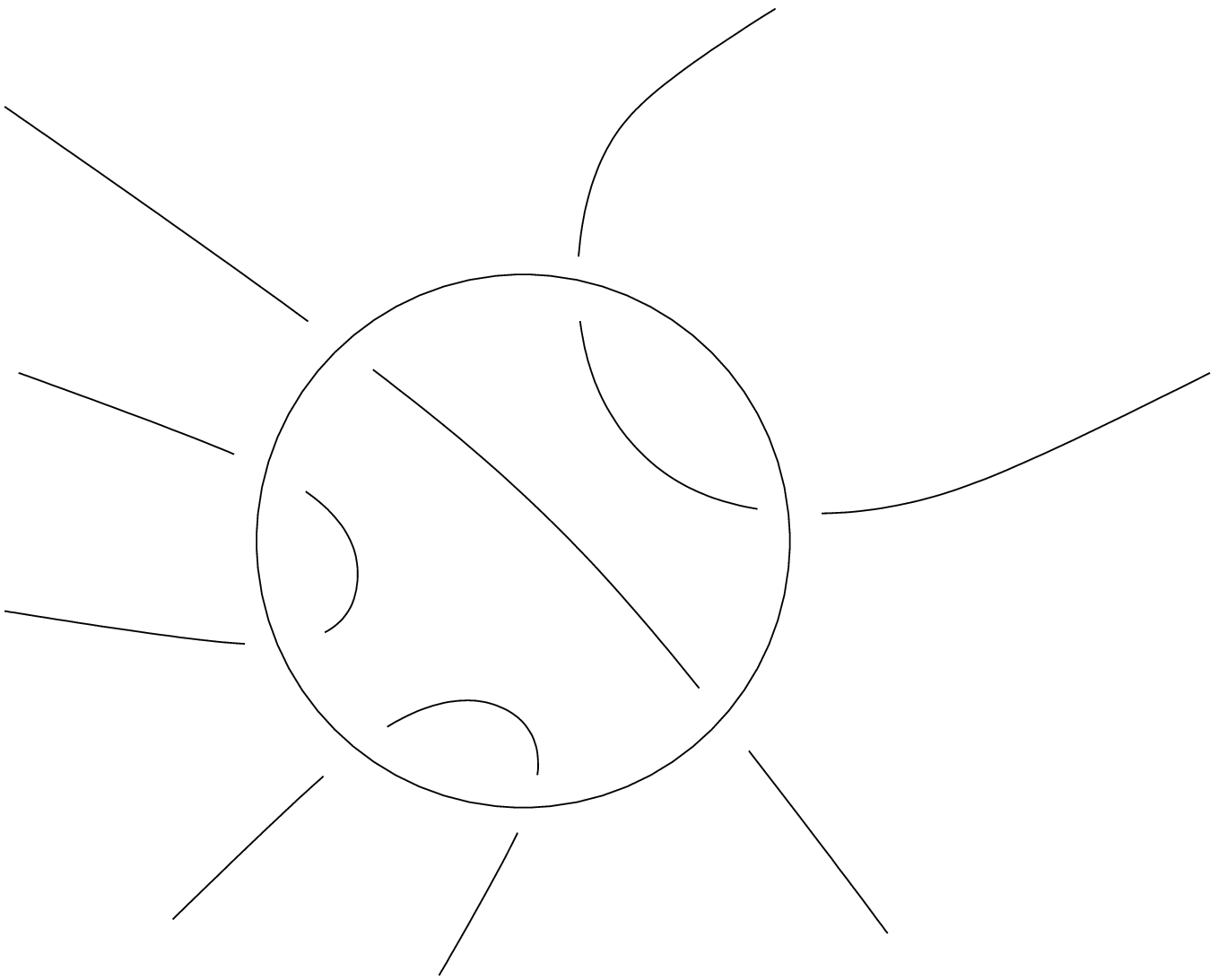}
\hskip 3truecm
\epsfxsize = 5truecm
\epsffile {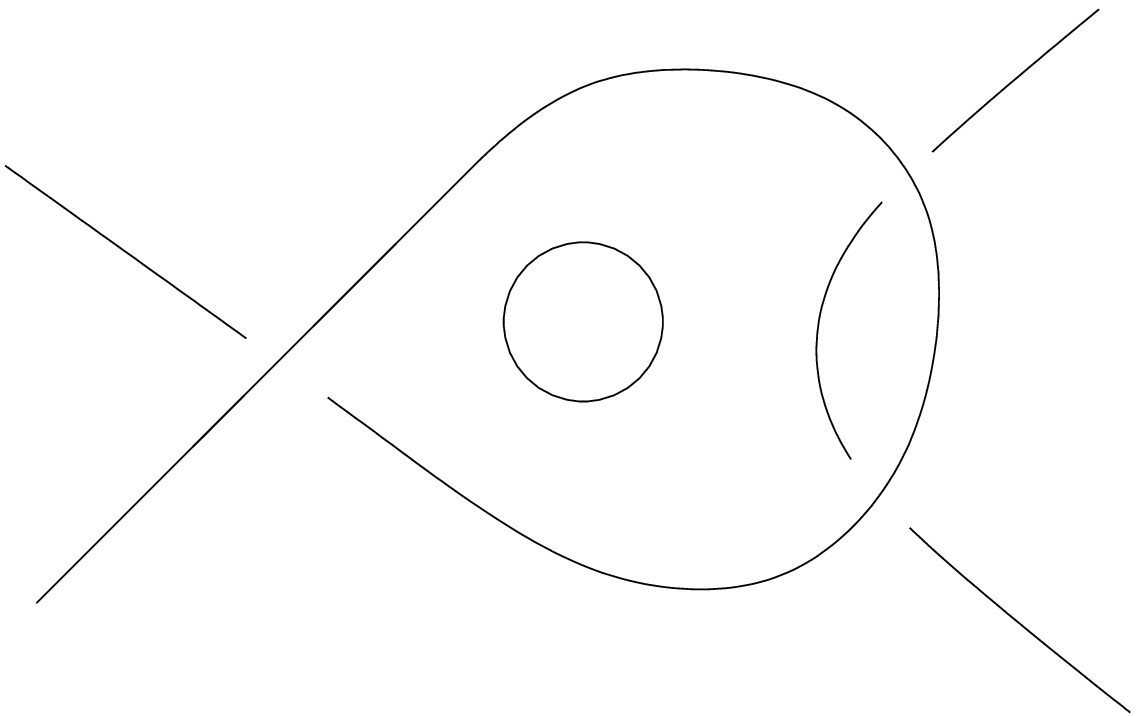}}}
\medskip
\hbox {\hskip 4truecm (a) \hskip 7.5truecm (b)}
\bigskip
\centerline {\hbox {
\epsfxsize = 5truecm
\epsffile {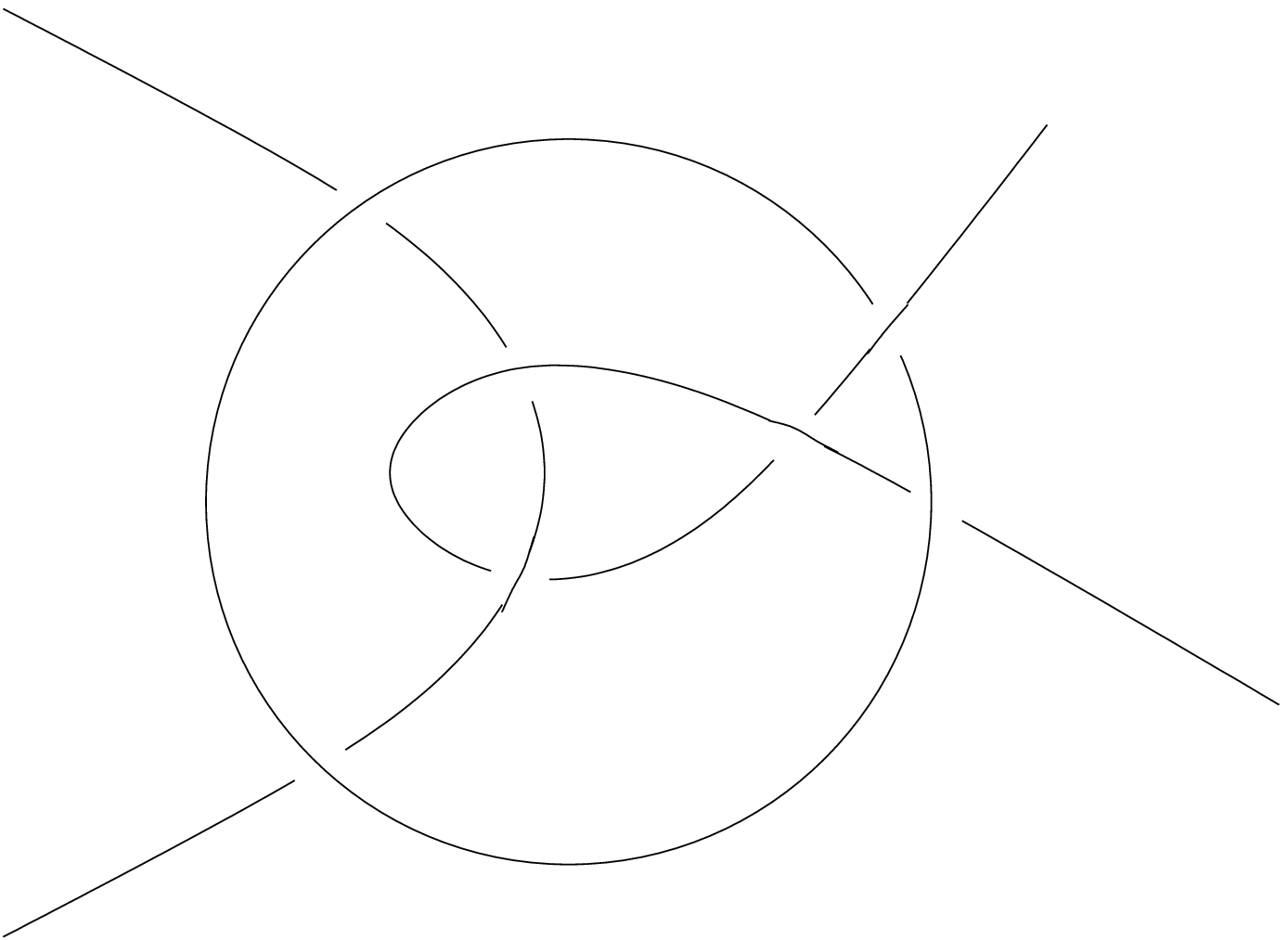} 
\hskip 3truecm
\epsfxsize = 2.5truecm
\epsffile {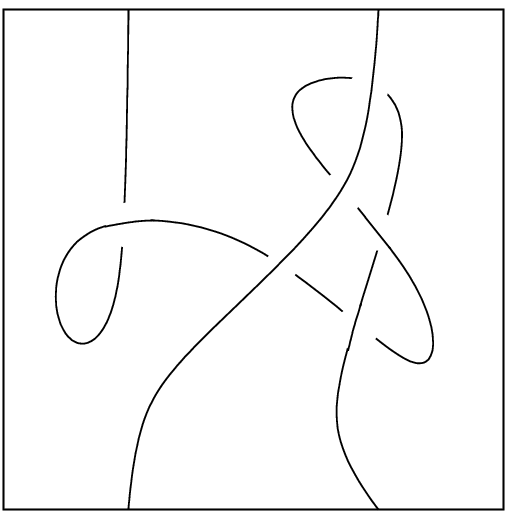}}}
\medskip
\hbox {\hskip 5truecm (c) \hskip 6.5truecm (d)}

\bigskip
\centerline {Figure 3}
\bigskip

\medskip
\noindent
{\bf Lemma 1.1.} \ If a given diagram of a tangle $T$ has $N$ crossings
and contains an improper bridge $b$ of length $B$, then the
diagram can be altered by type II and III Reidemeister moves
to a diagram with $N^\prime$ crossings and a bridge $b^\prime$
of length $B^\prime$ in such a way that $N^\prime<N$ and 
$N^\prime - B^\prime \le N - B$.

\medskip
\noindent 
{\it Proof:} \ For improper bridges in (a), (b), and (c) above, the
proof is identical to Lemma~0 in [2]. In the case of an
overcrossing arc, we may move the arc to the far right
of the tangle, leaving at most one crossing on that arc.
The two possible cases are shown in Figure~4.  The reduction
in crossing number is matched by a reduction in the length
of the bridge.
\square

\bigskip
\bigskip
\centerline {\hbox {
\epsfysize = 2.5truecm
\epsffile {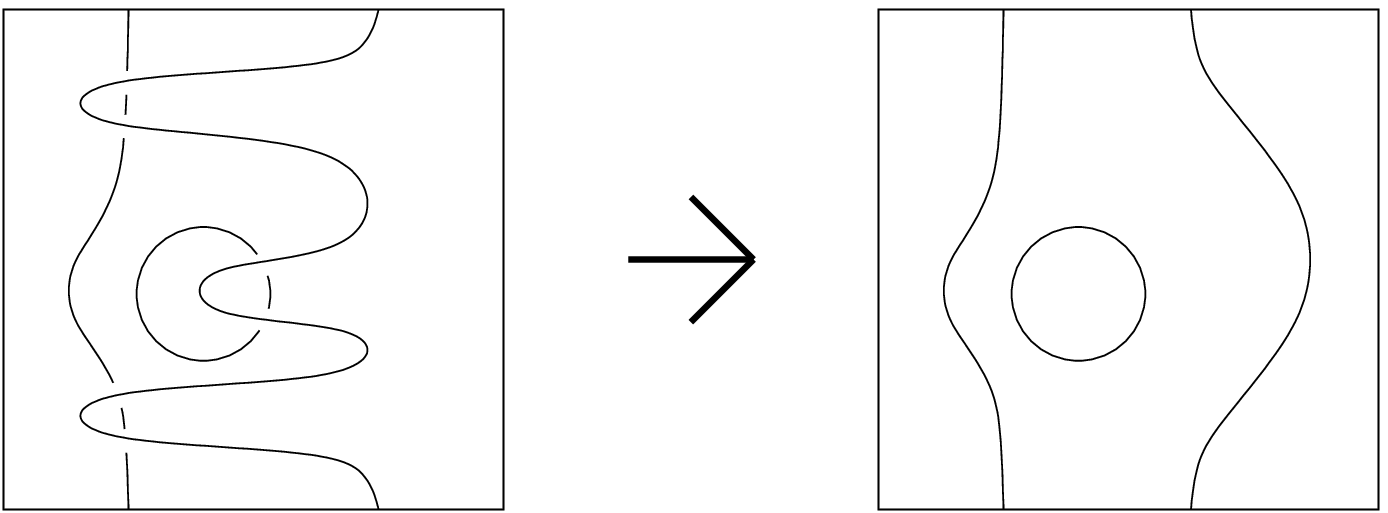}
\hskip 2truecm
\epsfysize = 2.5truecm
\epsffile {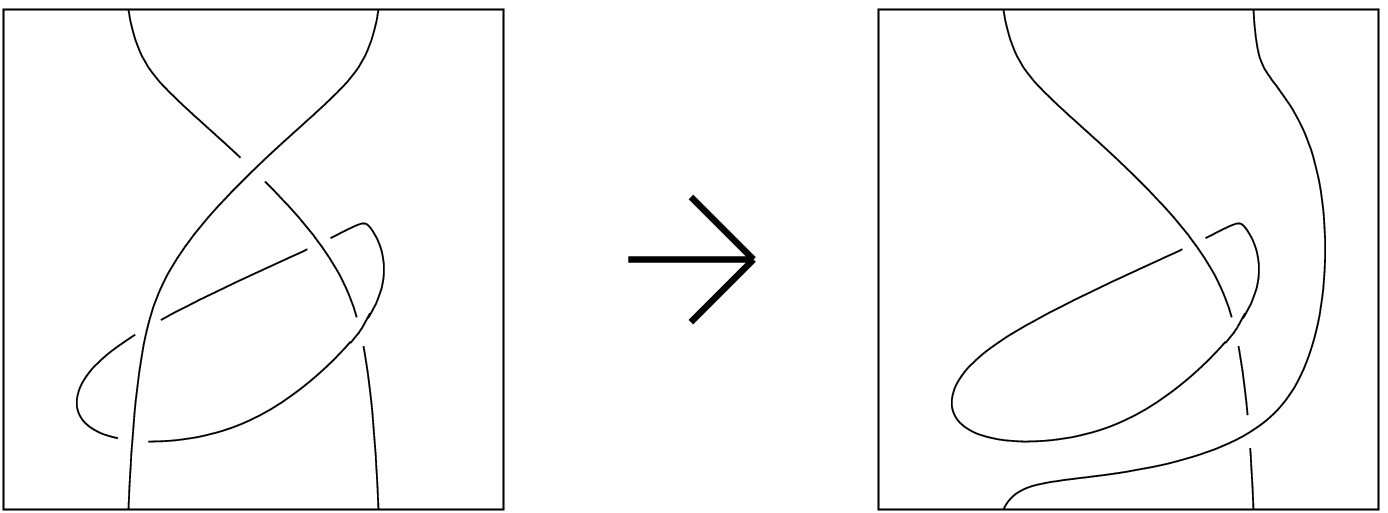}}}

\bigskip 
\centerline {Figure 4}
\bigskip

\bigskip
\noindent
{\bf Theorem 1.2.} \ Let $T$ be a tangle with $N$ crossings and a
bridge of length $B$.  Considered as an element of $M_2$,
$T$ can be written as a $\Lambda^\prime$-linear combination
$$f_1 (z, \lambda^{\pm 1}, \delta) P
+f_2 (z, \lambda^{\pm 1}, \delta) Q
+f_3 (z, \lambda^{\pm 1}, \delta) R_1
+f_4 (z, \lambda^{\pm 1}, \delta) R_2$$
where the $z$-degree of each $f_i$ is less than or
equal to $N-B$.

\medskip
\noindent
{\it Proof:} \ The theorem is true for any tangle with crossing
number $0$ or $1$ simply by replacing each disjoint circle
with a $\delta$ and a loop with $\lambda^{\pm 1}$.  If the
set of counterexamples to the theorem is nonempty, let $T$
be such a tangle diagram with minimal crossing number $N$;
among all such tangles with minimal crossing number, let
$T$ be one with maximal longest bridge length $B$.  Let $b$
be one of the longest bridges in $T$.  By Lemma 1.1, $b$
cannot be improper; thus it must have at least one endpoint
$c$ in the interior of $T$.  If we perform the skein
operation (i) at $c$, the two smoothings $T^0$ and $T^\infty$
will have smaller crossing number than $T$, while the new
$T^{\pm 1}$ will have a longer bridge.  None of these
tangles can be counterexamples to the theorem.  Thus the
$z$-degree of any coefficient 
of $T^0$ or $T^\infty$ is at most $(N-1)-B$, and
the $z$-degree of any coefficient of the changed
tangle $T^{\pm 1}$ is at most $N - (B+1)$. Thus there can
be no coefficient of $z T^\infty$ or $z T^0$ or
$T^{\pm 1}$ with $z$-degree high enough for $T$ to
be a counterexample to the theorem.
\square
\medskip

If one should want to write $T$ as a $\Lambda^\prime$-linear
combination of the basis tangles $P,Q,R_1$, one can do so
at the expense of adding $1$ to the $z$-degree of the
$P$ and $Q$ coefficients.  Note that $N-B=0$ in each of
the tangles $P,Q,R_1$, so performing a skein move to
replace $R_2$ with $z P - zQ + R_1$ adds $1$ to the $z$-degree
and accomplishes nothing as far as crossings and bridges 
are concerned.

We now consider a {\it wiring diagram} in the plane as
defined in Morton [3].  The endpoints of several disjoint
tangles are joined by non-crossing arcs.  See for example
Figure~5.  One obtains a link by inserting a tangle into
every box of a wiring diagram.

\bigskip
\centerline {
\epsfysize = 8truecm
\epsffile {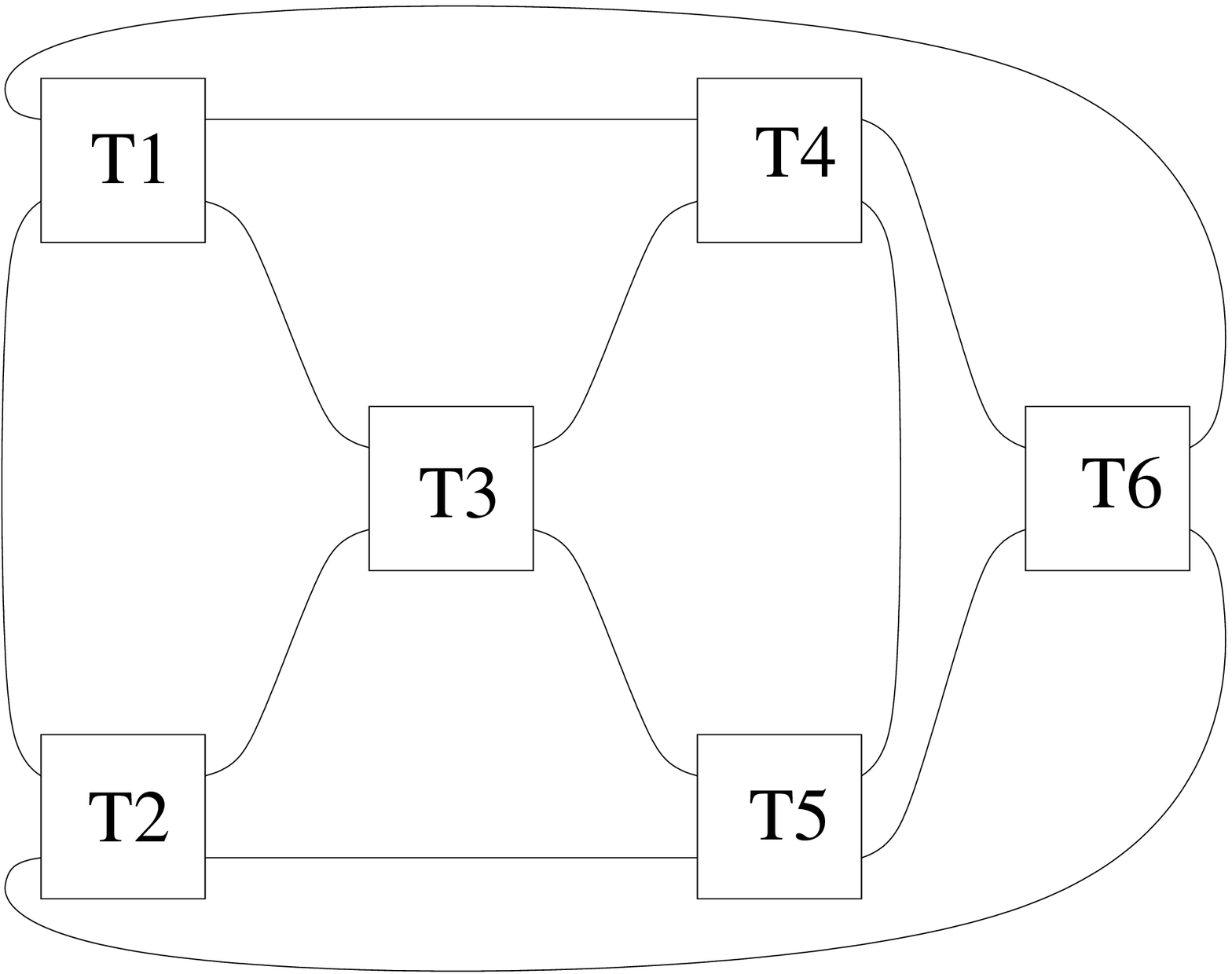}}

\bigskip 
\centerline {Figure 5}
\bigskip

\bigskip
\noindent
{\bf Theorem 1.3.} \ Let $L$ be a link diagram
written as a wiring diagram with $k$ tangles $\{T_i\}_1^k$.
Let tangle $T_i$ have $N_i$ crossings and a longest bridge
of length $B_i$.  Then the Dubrovnik polynomial
$D_L (\lambda^{\pm 1}, \delta, z)$ has $z$-degree less than
or equal to
$\displaystyle k-1+\sum_{i=1}^k (N_i - B_i)$

\medskip
\noindent
{\it Proof:} \ 
Apply the skein relations inside of each tangle $T_i$ to obtain
a linear combination of $P,Q,R_1$, and $R_2$ 
with $z$-degree bounded by
$N_i - B_i$. The result is that the Dubrovnik polynomial has
been written as a $\Lambda^\prime$-linear combination
of the polynomials of at most $4^k$ links, with each
coefficient of the combination having $z$-degree at most
$\sum_{i=1}^k (N_i - B_i)$.  Each of these $4^k$ links has
at most $k$ crossings and, unless it has no crossings at all,
a bridge of length at least $1$.  By the main theorem of [2],
each link has a Dubrovnik polynomial with $z$-degree at most
$k-1$.
\square
\medskip

It is customary to use relation (*) to replace the variable
$\delta$ in $\Lambda^\prime$ with 
$(\lambda^{-1} - \lambda)z^{-1}+1$.
Since this replacement introduces no positive powers of $z$ and
leaves one constant term, it has no effect on the
$z$-degree of a polynomial in $\Lambda^\prime$.  It is also
irrelevant to the $z$-degree of $D$, whether one uses the
regular isotopy invariant or the ambient isotopy invariant
version of the polynomial, since these differ only by a power
of $\lambda$.

\bigskip
\noindent
{\bf 2.  Examples from Rational Tangles.}
\bigskip

Consider the closed chain $L$
represented by the two different diagrams in Figure~6. The
seven-component link $L$
must have crossing number at least $14$, since each component
has linking number $\pm 1$ with exactly two other components
and so each component must have at least four crossings.
Therefore the diagrams in Figure~6 are minimal crossing
diagrams.  If we apply Theorem~1.3 to the top diagram,
we find two bridges of length $2$ which can be separated into
two separate tangles, and we therefore obtain an upper bound
of $14-2-1=11$ for the $z$-degree of this link.  However,
the bottom diagram may be decomposed into three tangles,
each containing a bridge of length $2$, and so we get an
upper bound of $14-3-1=10$.  The $z$-degree of $L$
is in fact $10$.

\bigskip
\centerline {
\epsfysize = 5truecm
\epsffile {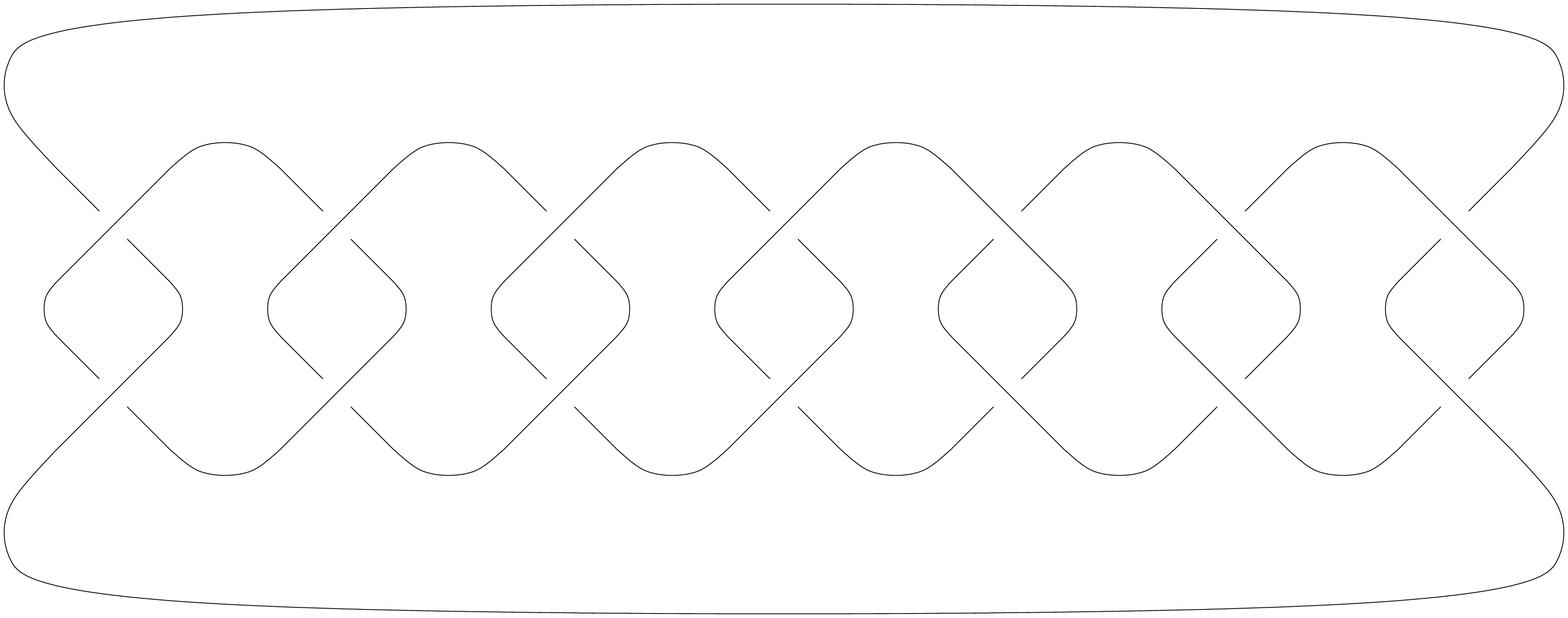}}

\bigskip
\centerline {
\epsfysize = 5truecm
\epsffile {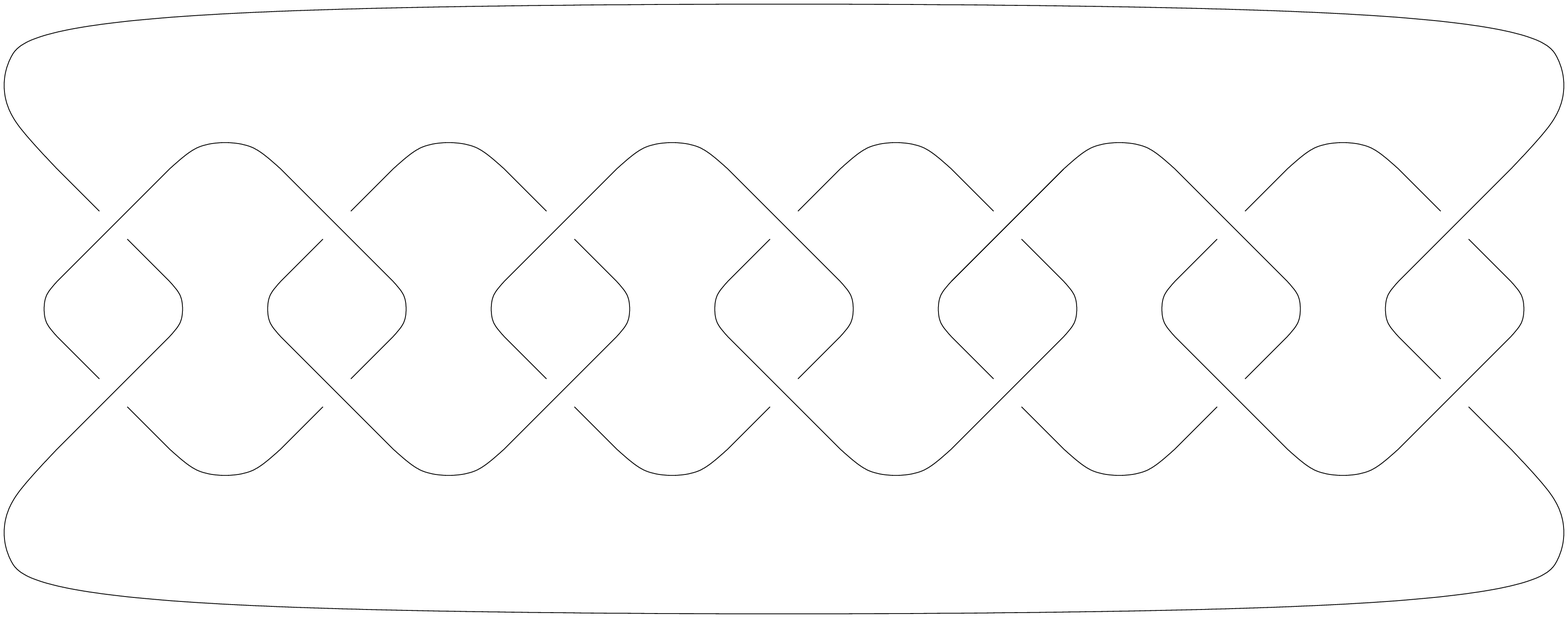}}

\bigskip 
\centerline {Figure 6}
\bigskip

More generally, let $L$ be a chain similar to the one in 
Figure~6
but with $p+q$ components.  Let $p$ be the number of positive
linkages (consecutive components ``positively'' linked)
and $q$ be the number of negative linkages.
For example, the top diagram of Figure~6 has four positive
linkages followed by three negative linkages, whereas the
bottom diagram alternates between positive and negative.
Arranging the linkages in positive-negative pairs, we
find that Theorem~1.3 gives a bound of $N - {\rm min\ }(p,q) -1$
for the $z$-degree of $L$ (where ${\rm min\ }(p,q)$ is
the smaller of the two numbers $p$ and $q$).
This bound is in fact exact.
However, grouping all the positive
linkages together and all the negative linkages together,
we get a bound of $N - 3$, which is far from exact if
${\rm min\ }(p,q)$ is large.

If $p=0$ or $q=0$, then $L$ is an alternating link,
and $N - {\rm min\ }(p,q)-1$ 
reduces to the $N-1$ of [2].
Consider the case $p>q>0$, so that
$N - {\rm min\ }(p,q) = N-q-1 = 2(p+q)-q-1 = 2p+q-1$.
Apply the skein relation (i)
to one of the crossings in a positive linkage, with
$L = L^+$.  Then $L^-$ is a connected sum of $p+q$
Hopf links, and so has degree $p+q-1 < 2p+q-1$.
$L^\infty$
is isotopic (picking up a factor of $\lambda$, of course)
to a chain with $N-2$ crossings, 
$p-1$ positive linkages, and $q$ negative
linkages, so it has $z$-degree $\le N - q -3$.  
Applying the move in Figure~7, $L^0$ is isotopic
(picking up $\lambda^{-1}$) to a chain with
$p$ positive linkages and $q-1$ negative linkages.
Inductively, we see that the $L^0$ term has the 
largest $z$-degree, and that the $z$-degree of
$L$ is therefore $N-q-1$.  

\bigskip
\centerline {
\epsfysize = 4truecm
\epsffile {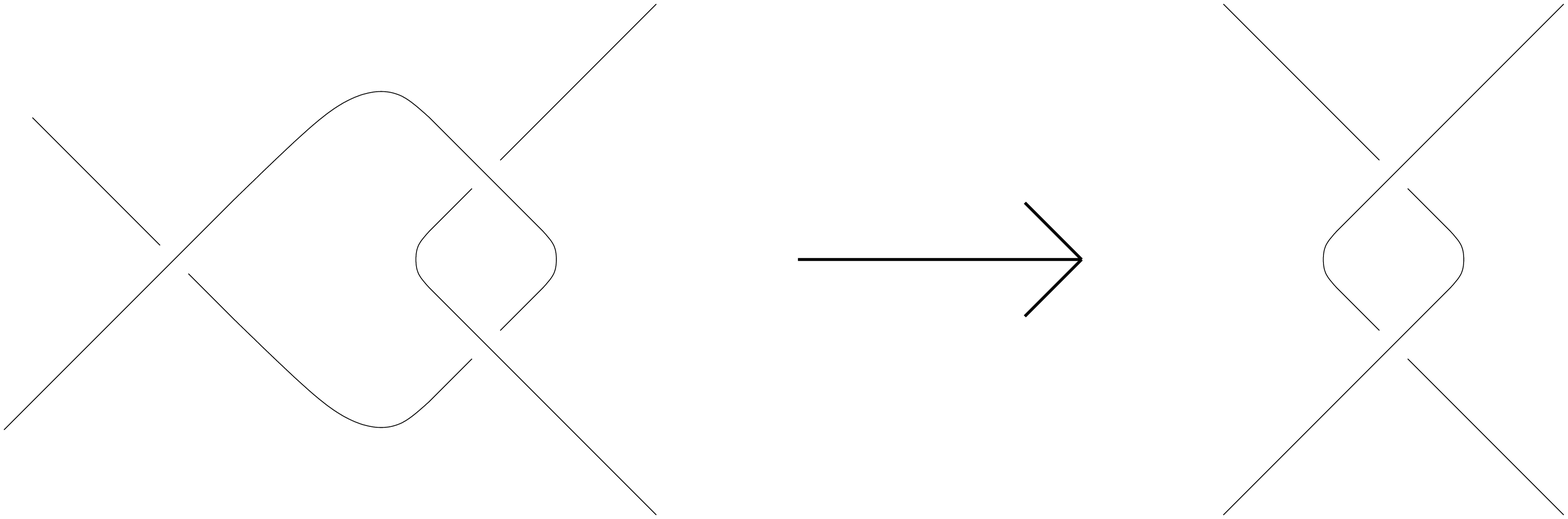}}

\bigskip 
\centerline {Figure 7}
\bigskip

The case $p<q$ follows similarly.  The case $p=q$ 
is more complicated to do by this method because
in applying the skein relation (i), one encounters
the case that more than one of the three replacement
terms have maximal degree.  We will see that the
formula holds in the $p=q$ case below.  
First, however, we want to point out that the chains
described so far are all prime links, so that the
reduction in $z$-degree obtained (reduction below the
crossing number) is not just the connected sum
effect (as mentioned in the introduction) in disguise.
For suppose that such a link $L$ has a two-sphere $S$
which intersects $L$ in exactly two points.  Those
two points must be on the same component $K$.  
$L-K$ cannot be split, no matter what $K$ is
(linking number considerations again), so every
component of $L-K$ lies on the same side of $S$.
On the other side of $S$ must be
only a single arc of $K$, which cannot
be knotted because none of the components of $L$ are
individually knotted.

More examples may be obtained by considering chains
where two consecutive components may be twisted together
any even number of times.  Such a chain with $k$
components may be indexed by $k$ nonzero even integers
$(m_1, m_2, \dots m_k)$.  The integers describe
the linking between consecutive components.  For
example, the diagram indexed by $(2,4,-4,2)$ is
shown in Figure~8.  The same arguments used above
show that such a diagram has minimal crossing number
and represents a prime link.  If the $m_i$ are permuted,
then even if the represented link changes, the Kauffman
polynomial does not, since permutations may be
accomplished by mutations.  Let $p$ be the number of
positive $m_i$ (which correspond to positive linkages)
and let $q$ be the number of negative
$m_i$ (which correspond to negative linkages),
so that $k = p+q$.
If the positive and negative linkages are arranged
in alternating fashion, then once again Theorem~1.3
gives us a bound of $N - {\rm min\ }(p,q) - 1$.  Again,
we shall see that this bound is exact.  Arranging
the positive linkages all together, however, again
gives a bound of $N - 3$.

\bigskip
\centerline {
\epsfysize = 7truecm
\epsffile {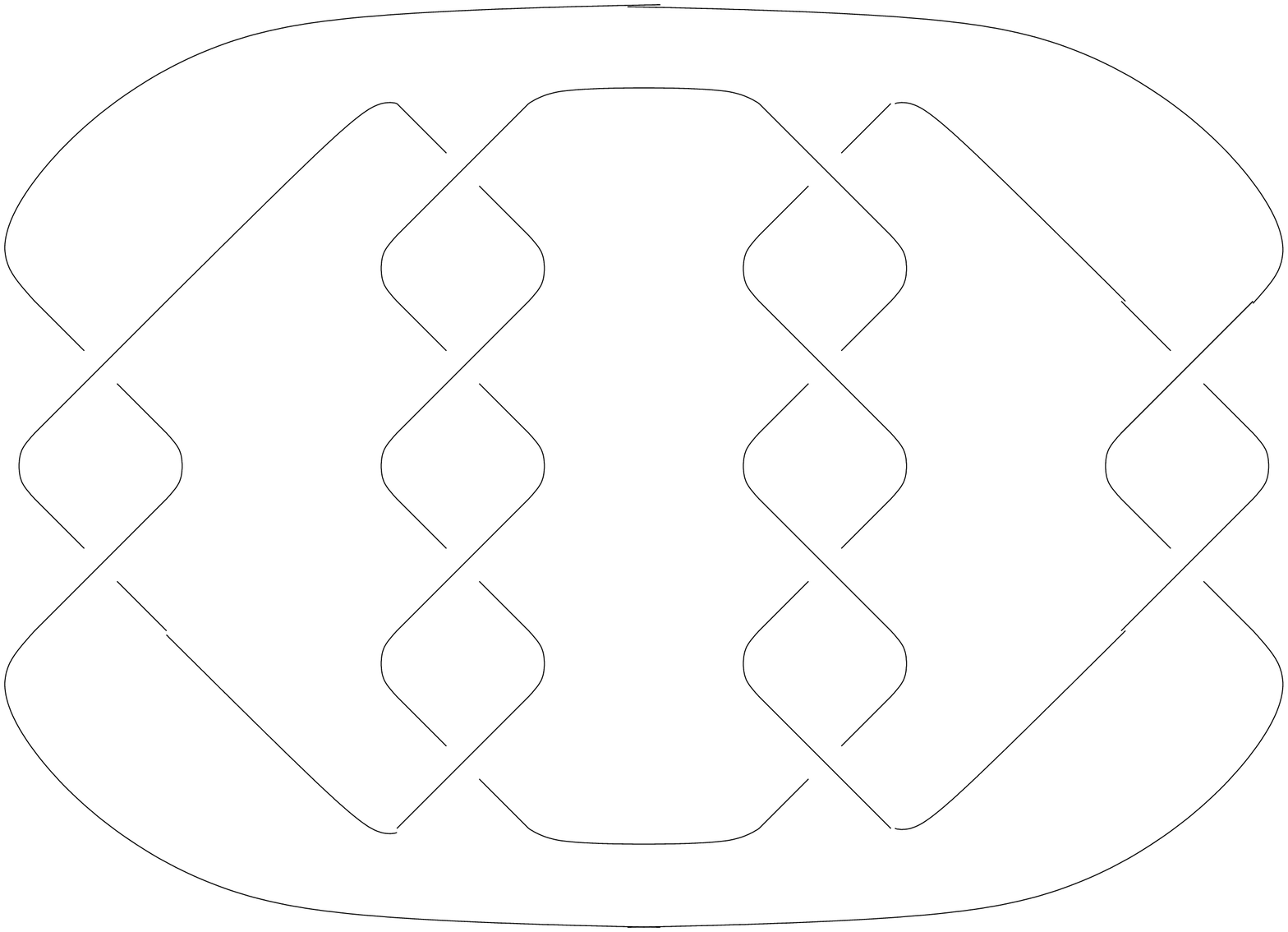}}

\bigskip 
\centerline {Figure 8}
\bigskip

We now define $R_1$ to be 
a {\it positive rational tangle}.  Moreover, 
we will say that any tangle built up out of a positive
rational tangle $T$ and $R_1$ in either of the two ways shown
in Figure~9 is also a positive rational tangle.  We
define a {\it negative rational tangle} similarly, building
with $R_2$ instead of $R_1$.  The term {\it rational}
comes from Conway [1], and in fact all rational
tangles are either positive or negative (in our sense),
depending on whether their associated continued fractions
(in Conway's sense) represent positive or negative numbers.
It is clear that a positive or negative rational tangle
is alternating.
A positive tangle that can be written
as in Figure~9a will be called {\it vertical}, 
and one that can
be written as in Figure~9b will be called {\it horizontal}.
We define $R_1$ to be neither horizontal nor vertical.
We define horizontal and vertical negative tangles similarly.

\bigskip
\centerline {
\epsfysize = 5truecm
\epsffile {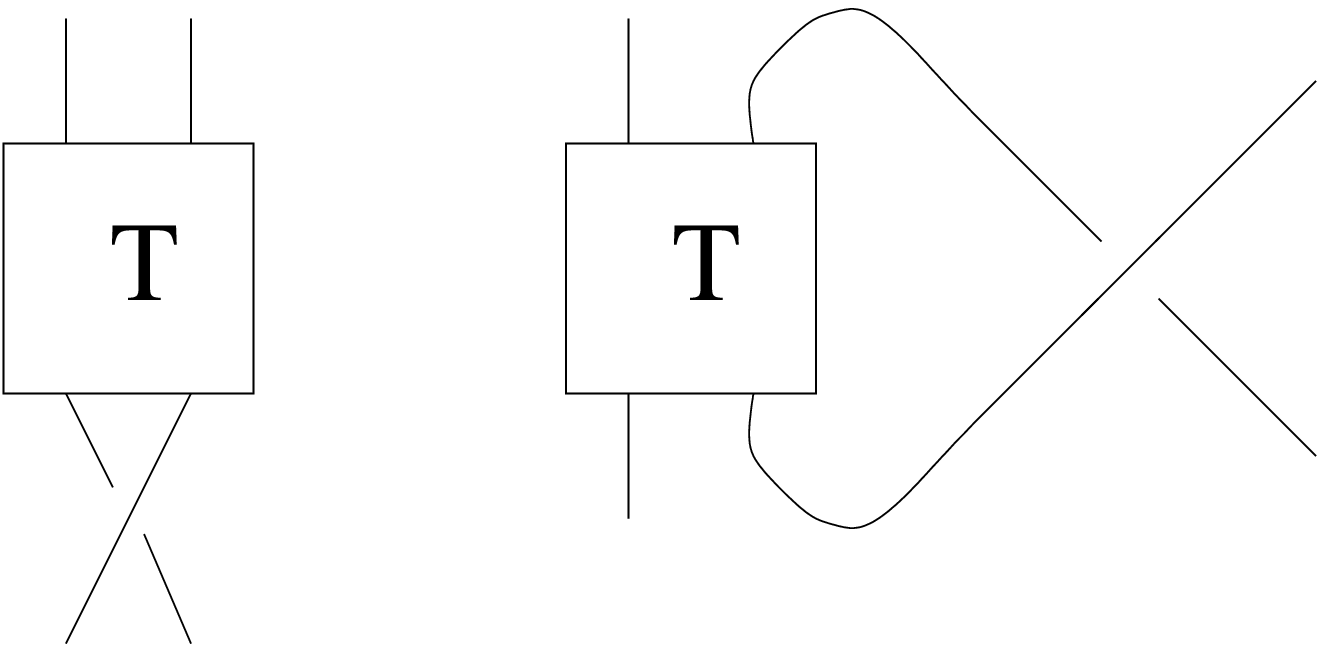}}
\medskip
\hbox {\hskip 3truecm (a) vertical \hskip 3truecm (b) horizontal}

\bigskip 
\centerline {Figure 9}
\bigskip

Another way of stating Theorem~1.2 is that any tangle may
be written as a polynomial in $z$ of
degree $\le N-B$, where
the coefficients of the polynomial are themselves polynomials
in $\lambda^{\pm 1}$, $P$, $Q$, $R_1$, and $R_2$.  For rational
tangles, we have the following:

\bigskip
\noindent
{\bf Theorem 2.1.}
A positive rational tangle with $N$ crossings
may be uniquely written as a polynomial
in $z$ of degree $N-1$, where the coefficients are polynomials
in $\lambda^{\pm 1}$, $P$, $Q$, and $R_1$.  If $N>1$, then
the $z^{N-1}$ coefficient of a vertical tangle is 
$\pm (R_1 -\lambda Q)$, and the $z^{N-1}$ coefficient of a
horizontal tangle is $\pm (R_1 - \lambda^{-1} P)$.
A negative rational tangle may be uniquely written as a polynomial
in $z$ of degree $N-1$, where the coefficients are polynomials
in $\lambda^{\pm 1}$, $P$, $Q$, and $R_2$.  If $N>1$, then
the $z^{N-1}$ coefficient of a vertical tangle is 
$\pm (R_2 -\lambda^{-1} Q)$, and the $z^{N-1}$ coefficient of a
horizontal tangle is $\pm (R_2 - \lambda P)$.

\medskip
\noindent
{\it Proof:}
The uniqueness follows from the result of Morton and Trazcyk [4],
that $P$, $Q$, and $R_1$ form a free basis for $M_2$.  (It is
clear that $R_1$ may be replaced by $R_2$ without affecting
this statement.)  Beyond that, it is a simple induction argument,
applying the skein relation (i) to $R_1$ in Figure~9
for the positive case (and similarly for the negative case).
\square
\medskip

Now it follows easily
that if $L$ is constructed by wiring $p$ positive
and $q$ negative vertical rational 
tangles together in a horizontal line,
the $z$-degree of $L$ is $N - {\rm min\ }(p,q) -1$.
We replace each rational tangle with a linear combination
as in the proof of 
Theorem~1.3.  If $p>q$, then the largest $z$-degree
in the result is obtained from the $R_1$ term in each
positive tangle and the $Q$ term in each negative tangle.
The leading $z$-coefficient will be $\pm \lambda^{-q}$
times the leading $z$-coefficient of the $(2,p)$ torus
link (which is built up of $p$ copies of $R_1$ and
$q$ copies of $Q$).
If $q<p$, then the largest $z$-degree is obtained from
the $R_2$ term in each negative tangle and the $Q$ term
in each positive tangle.  
The leading $z$-coefficient will be $\pm \lambda^{p}$
times the leading $z$-coefficient of the $(2,p)$ torus link.
If $p=q$, then both of the terms just described are
of maximal degree, but because they differ by 
$\lambda^{2p}$ they cannot cancel.

\bigskip
\noindent
{\bf References.}
\medskip

\smallskip
\item {[1]}
J.H. Conway, 
{\it An enumeration of knots and links, and some
of their algebraic properties}, 1970 Computational
Problems in Abstract Algebra (Proc. Conf., Oxford, 1967),
329--358.  Pergamon, Oxford.

\smallskip
\item {[2]}
M.E. Kidwell,
{\it On the degree of the Brandt-Lickorish-Millett-Ho
polynomial of a link}, Proceedings of the American Mathematical
Society 100 (1987), 755--762.

\smallskip
\item{[3]}
H.R. Morton,
{\it Invariants of links and 3-manifolds from skein theory
and from quantum groups}, Topics in Knot Theory (Erzurum, 1992),
107--155. NATO ASI Series C, 399, Kluwer Academic Publishers, 1993.

\smallskip
\item {[4]}
H.R. Morton and P. Traczyk,
{\it Knots and algebras},
Contribuciones Matematicas en homenaje al professor
D.~Antonio Plans Sanz de~Bremond (ed. E.~Martin-Peinador
and A.~Rodez Usan), University of Zaragoza (1990) 201--220.

\smallskip
\item {[5]}
D. Rolfsen,
{\it Knots and Links},
Mathematics Lecture Series, volume~7.
Publish or Perish, Inc., Wilmington, DE, 1976.

\smallskip
\item{[6]}
M.B. Thistlethwaite, 
{\it Kauffman's polynomial and alternating links}, \hfil \break
Topology 27 (1988), 311--318.

\smallskip
\item {[7]}
Y. Yokota,
{\it The Kauffman polynomial of alternating links},
Topology and its Applications 65 (1995), 229--236.

\end